\RequirePackage{ifpdf}
\ifpdf 
\documentclass[pdftex]{sigma}
\else
\documentclass{sigma}
\fi

\usepackage[all]{xy}
\begin{document}

\newcommand\ugarr{\!\!\!&=&\!\!\!}

\newcommand{\Mreg}{M_{\mathrm{reg}}}

\newcommand{\for}[1]{(\ref{#1})}

\newcommand\aHp{almost-Hamiltonian}
\newcommand\asp{almost-symplectic}
\newcommand\aPp{almost-Poisson}
\newcommand\aPBp{almost-Poisson brackets}
\newcommand\aH{almost-Hamiltonian\ }
\newcommand\as{almost-symplectic\ }
\newcommand\aPB{almost-Poisson brackets\ }
\newcommand\aP{almost-Poisson\ }

\newcommand{\bib}[1]{\textup{\cite{#1}}}

\newcommand{\fcomment}[1]{\marginpar{\scriptsize\tt #1 (F)}}

\def\commento#1{\marginpar{\tiny #1}}

\def\bull{\item[$\bullet$]}
\def\boxed#1{\fbox{$\displaystyle#1$}}
\def\boxedtext#1{\fbox{#1}}
\def\commento#1{\marginpar{\tiny #1}}
\newcommand{\vect}[1]{#1}   
\newcommand{\mezzo}{{1\over2}\, }

\newcommand{\step}[1]{^{(#1)}}

\renewcommand{\a}{\alfa}
\renewcommand{\b}{\beta}
\newcommand{\e}{\epsilon}
\newcommand{\g}{\gamma}
\newcommand{\G}{\Gamma}
\newcommand{\p}{\varphi}
\renewcommand{\P}{\Phi}
\newcommand{\s}{\sigma}
\renewcommand{\a}{\alpha}
\renewcommand{\d}{\delta}
\renewcommand{\o}{\omega}
\renewcommand{\r}{\rho}
\renewcommand{\l}{\lambda}
\renewcommand{\L}{\Lambda}
\renewcommand{\O}{\Omega}
\renewcommand{\s}{\sigma}    \renewcommand{\S}{\Sigma}


\newcommand{\beq}[1]{\begin{equation}\label{#1}}
\newcommand{\eeq}{\end{equation}}

\newcommand{\bWS}{\begin{WriteSmall}\acapo}
\newcommand{\eWS}{\acapo\end{WriteSmall}}

\newcommand{\Scrivipiccolo}[1]{{\oddsidemargin=2cm\footnotesize #1} }

\newcommand{\cA}{\mathcal{A}}
\newcommand{\cB}{\mathcal{B}}
\newcommand{\cC}{\mathcal{C}}
\newcommand{\cD}{\mathcal{D}}
\newcommand{\cE}{\mathcal{E}}
\newcommand{\cF}{\mathcal{F}}
\newcommand{\cG}{\mathcal{G}}
\newcommand{\cH}{\mathcal{H}}
\newcommand{\cK}{\mathcal{K}}
\newcommand{\cI}{\mathcal{I}}
\newcommand{\cL}{\mathcal{L}}
\newcommand{\cM}{\mathcal{M}}
\newcommand{\cO}{\mathcal{O}}
\newcommand{\cP}{\mathcal{P}}
\newcommand{\cQ}{\mathcal{Q}}
\newcommand{\cR}{\mathcal{R}}
\newcommand{\cS}{\mathcal{S}}
\newcommand{\cT}{\mathcal{T}}
\newcommand{\cU}{\mathcal{U}}
\newcommand{\cV}{\mathcal{V}}
\newcommand{\cW}{\mathcal{W}}
\newcommand{\cX}{\mathcal{X}}

\newcommand{\deq}{$\interi_4\per\interi_2$--equivalent\ }
\newcommand{\deqp}{$\interi_4\per\interi_2$--equivalent}
\newcommand{\deqence}{$\interi_4\per\interi_2$--equivalence\ }
\newcommand{\deqencep}{$\interi_4\per\interi_2$--equivalence}
\newcommand{\naturali}{\mathbb{N}}
\newcommand{\Sp}[1]{ \mathrm{Sp}(#1)}
\newcommand{\fou}[2]{\langle #1\rangle_{#2}}


\newcommand{\acapo}{\vspace{2ex}}
\newcommand{\acapon}{\vspace{2ex}\noindent}
\newcommand{\qqquad}{\quad\quad\quad}

\newcommand{\ins}[3]{\vbox to0pt{\kern-#2 \hbox{\kern#1
#3}\vss}\nointerlineskip}

\newcommand{\toro}{\mathbb{T}}
\newcommand{\complessi}{\mathbb{C}}
\newcommand{\reali}{\mathbb{R}}
\newcommand{\razionali}{\mathbb{Q}}
\newcommand{\interi}{\mathbb{Z}}

\newcommand{\tdue}{\toro^2}
\newcommand{\tenne}{\toro^n}
\newcommand{\rdue}{\reali^2}
\newcommand{\rtre}{\reali^3}
\newcommand{\rsei}{\reali^6}
\newcommand{\renne}{\reali^n}

\newcommand{\unity}{\mathbf{1}}


\newcommand{\pv}{\per}
\newcommand{\ug}{\;=\;}
\newcommand{\piu}{\;+\;}
\newcommand\meno{\;-\;}
\newcommand\TO{\;\to\;}
\newcommand\Mapsto{\;\mapsto\;}
\newcommand\GE{\;\ge\;}
\newcommand\LE{\;\le\;}

\newcommand\IMPLY{\quad\Longrightarrow\quad}

\newcommand{\per}{\times}
\renewcommand{\bar}{\overline}

\newcommand{\der}[2]{\frac{\partial#1}{\partial#2}}
\newcommand{\dder}[3]{\frac{\partial^2#1}{\partial#2\partial#3}}

\newcommand{\grad}{\mathrm{grad}}
\newcommand{\rot}{\mathrm{rot}}
\newcommand{\cost}{\mathrm{cost}}
\newcommand{\const}{\mathrm{cost}}
\newcommand\rank{\mathrm{rank\,}}
\newcommand\diag{\mathrm{diag\,}}
\newcommand\tr{\mathrm{tr\,}}
\renewcommand\skew{\mathrm{skew\,}}
\newcommand\sign{\mathrm{sign\,}}

\newcommand\fineprova{\hfill \vrule height6pt width6pt depth0pt}

\newcommand{\SKEWTRE}{\mathrm{so(3)}}
\newcommand{\OTRE}{\mathrm{O(3)}}
\newcommand{\LTRE}{\mathrm{L(3)}}
\newcommand{\GLTRE}{\mathrm{GL(3)}}
\newcommand{\SLTRE}{\mathrm{SL(3)}}
\newcommand{\sltre}{\mathrm{sl(3)}}
\newcommand{\SLDTRE}{\mathrm{SLD(3)}}
\newcommand{\sldtre}{\mathrm{sld(3)}}
\newcommand{\SOTRE}{\mathrm{SO(3)}}
\newcommand{\sotre}{\mathrm{so(3)}}
\newcommand{\PS}[2]{\langle#1,#2\rangle}
\newcommand{\bigPS}[2]{\big\langle#1,#2\big\rangle}
\newcommand{\BigPS}[2]{\Big\langle#1,#2\Big\rangle}
\newcommand{\colvect}[2]{\left(\matrix{#1\cr#2\cr}\right)}

\newcommand{\PP}[2]{\{#1,#2\}}
\newcommand{\PPP}[3]{\PP{#1}{\PP{#2}{#3}}}

\newcommand{\sub}[1]{_{(#1)}}


\newcommand{\D}{D}
\newcommand{\MR}{G^R}      
\newcommand{\MS}{G^S}      
\newcommand{\wMS}{G}       
\newcommand{\wMR}{G}       
\newcommand{\GG}{\Psi}     

\newcommand{\Sdue}{$\mathrm{S}_2$ }
\newcommand{\Stre}{$\mathrm{S}_3$ }
\newcommand{\Sduep}{$\mathrm{S}_2$}
\newcommand{\Strep}{$\mathrm{S}_3$}

\newcommand{\cBSdue}{\cB_{\mathrm{S}_2}}
\newcommand{\cBStre}{\cB_{\mathrm{S}_3}}
\newcommand{\cBI}{\cB_{\mathrm{I}}}
\newcommand{\cBII}{\cB_{\mathrm{II}}}
\newcommand{\cBIII}{\cB_{\mathrm{III}}}

\newcommand{\twopig}{}


\newcommand{\mytable}[4]{
\begin{table}
\center{\small
\begin{tabular}{||l|l||l|l||l|l||l|l||l|l||}
\hline\hline
\multicolumn{10}{||c||}{}\\
\multicolumn{10}{||c||}{{\em #2 }} \\
\multicolumn{10}{||c||}{}\\
\multicolumn{1}{||c|}{#3} & \multicolumn{1}{|c||}{#4} &
\multicolumn{1}{c|}{#3} & \multicolumn{1}{|c||}{#4} &
\multicolumn{1}{c|}{#3} & \multicolumn{1}{|c||}{#4} &
\multicolumn{1}{c|}{#3} & \multicolumn{1}{|c||}{#4} &
\multicolumn{1}{c|}{#3} & \multicolumn{1}{|c||}{#4} \\
\hline
\input{#1}
\hline\hline
\end{tabular}
}
\end{table}
}

\newcommand{\E}{{\mathbb E}}
\newcommand{\F}{{\mathbb F}}
\newcommand{\K}{{\mathbb K}}
\newcommand{\J}{{\mathbb J}}

\def\mommap{\mathrm{\bf J}}


\newcommand{\norm}[1]{\left \vert #1 \right \vert}
\newcommand{\proj}[2]{{\mathbb P}_{#1} #2}
\newcommand{\la}{\left \langle}
\newcommand{\ra}{\right \rangle}
\newcommand{\lsb}{\left [}
\newcommand{\rsb}{\right ]}
\newcommand{\lcb}{\left \{}
\newcommand{\rcb}{\right \}}

\newcommand{\sands}{\qquad \mbox{and} \qquad}
\newcommand{\dcomment}[1]{\footnote{{\sl #1 (D)}}}
\newcommand{\half}{{\textstyle {1 \over 2}}}
\newcommand{\smallfrac}[2]{{\textstyle {#1 \over #2}}}

\newcommand{\ddt}{\smallfrac{d\ }{dt}}
\newcommand{\dds}{\smallfrac{d\ }{ds}}
\newcommand{\ddeps}{\smallfrac{d\ }{d\e}}
\newcommand{\ddlambda}{\smallfrac{d\ }{d\l}}


\newcommand\bList{
\begin{list}{}{\leftmargin2em\labelwidth1em\labelsep.5em\itemindent0em
\topsep0ex\itemsep-.8ex} }
\newcommand\eList{\end{list}}

\newcommand\bListList{
\begin{list}{}{\leftmargin2.5em\labelwidth1.7em\labelsep.5em\itemindent0em
\topsep0ex\itemsep-.8ex} }
\newcommand\eListList{\end{list}}

\newcommand\finelemma{\raisebox{4pt}
           {\framebox[6pt][5pt]{\hbox to 6pt{}}
}}

\def\uncatcodespecials{\def\do##1{\catcode`##1=12 }\dospecials}
\def\sic{\begingroup\tt\uncatcodespecials\obeylines
  \obeyspaces\doverbatim}
\newcount\balance
{\catcode`<=1 \catcode`>=2 \catcode`\{=12 \catcode`\}=12
  \long\gdef\doverbatim{<\balance=1\verbatimloop>
  \long\gdef\verbatimloop#1<\def\next<#1\verbatimloop>%
    \if#1{\advance\balance by1
    \else\if#1}\advance\balance by-1
     \ifnum\balance=0\let\next=\endgroup\fi\fi\fi\next>>
%


\def\de#1{\partial_{#1}}
\def\dede#1#2{\frac{\partial #1}{\partial #2}}


\newcommand{\arcsech}{{\rm arcsech}}
\renewcommand{\choose}[2]{\lp\genfrac{}{}{0pt}{1}{#1}{#2}\rp}
\newcommand{\cis}{{\rm cis}\ \!}


\newcommand{\her}{\rm Her}
\newcommand{\sher}{\rm SkewHer}
\newcommand{\antiher}{\rm AntiHer}
\newcommand{\symm}{{\rm symm}}
\newcommand{\reg}{{\rm reg}}
\newcommand{\Cart}{\rm Cart}
\newcommand{\SO}{{\rm SO}}
\newcommand{\SL}{{\rm SL}}
\newcommand{\SU}{\rm SU}
\newcommand{\LieAlg}{\mathfrak g}
\renewcommand{\c}{\mathfrak c}
\newcommand{\CartanAlg}{\mathfrak t}
\renewcommand{\s}{\mathfrak s}
\renewcommand{\k}{\mathfrak k}
\newcommand{\z}{\mathfrak z}
\newcommand{\h}{\mathfrak h}
\newcommand{\orth}{\mathfrak o}
\newcommand{\so}{\mathfrak so}
\newcommand{\uni}{\mathfrak u}
\newcommand{\suni}{\mathfrak su}

\newcommand{\im}{{\ \!\rm im}\ \!}
\renewcommand{\diag}{{\rm diag}\ \!}
\newcommand{\ad}{{\rm ad}}
\newcommand{\Ad}{{\rm Ad}}
\renewcommand{\tr}{{\rm tr}\ }
\newcommand{\eval}{{\rm e-val}\ }
\newcommand{\evect}{{\rm e-vect}\ }
\newcommand{\orbit}{{\mathcal O}}
\newcommand{\emmap}{{\mathcal EM}}
\newcommand{\emdomain}{{\mathcal D}}

\renewcommand{\rank}{{\rm rank}\ \!}
\newcommand{\codim}{{\rm codim}\ \!}


\renewcommand{\l}{\lambda}
\newcommand{\ve}{\varepsilon}
\newcommand{\vp}{\varphi}
\newcommand{\vt}{\vartheta}
\renewcommand{\const}{{\rm const}}
\newcommand{\id}{{\rm id}}
\def\ml{l\kern-0.02cm$^\prime$\kern-0.03cm}


\newcommand{\open}[1]{\overset\circ{#1}}
\newcommand{\ul}{\underline}
\newcommand{\wt}{\widetilde}
\newcommand{\wh}{\widehat}
\newcommand{\pr}{{\rm pr}\ \!}
\newcommand{\st}{\ |\ }


\newcommand{\R}{{\mathbb R}}
\newcommand{\Z}{{\mathbb Z}}
\newcommand{\C}{{\mathbb C}}
\newcommand{\Q}{{\mathbb Q}}
\newcommand{\T}{{\mathbb T}}
\newcommand{\I}{{\mathbb I}}


\renewcommand{\[}{\left[}
\renewcommand{\]}{\right]}
\newcommand{\lp}{\left(}
\newcommand{\rp}{\right)}

\def\matrice#1{\left(\begin{smallmatrix}#1\end{smallmatrix}\right)}

\allowdisplaybreaks

\renewcommand{\PaperNumber}{051}

\renewcommand{\thefootnote}{$\star$}

\FirstPageHeading

\ShortArticleName{Geometry of Certain Integrable Systems with
Symmetry}

\ArticleName{Geometry of Invariant Tori of Certain Integrable\\
Systems with Symmetry and an Application\\ to a Nonholonomic
System\footnote{This paper is a contribution to the Proceedings of
the Workshop on Geometric Aspects of Integ\-rable Systems
 (July 17--19, 2006, University of Coimbra, Portugal).
The full collection is available at
\href{http://www.emis.de/journals/SIGMA/Coimbra2006.html}{http://www.emis.de/journals/SIGMA/Coimbra2006.html}}}

\Author{Francesco FASS\`O and Andrea GIACOBBE}
\AuthorNameForHeading{F. Fass\`o and A. Giacobbe}

\Address{Dipartimento di Matematica Pura e Applicata, Universit\`a di Padova,\\
Via Trieste 63, 35131 Padova, Italy}
\Email{\href{mailto:fasso@math.unipd.it}{fasso@math.unipd.it},
\href{mailto:giacobbe@math.unipd.it}{giacobbe@math.unipd.it}}

\ArticleDates{Received November 20, 2006, in f\/inal form March
15, 2007; Published online March 22, 2007}

\Abstract{Bif\/ibrations, in symplectic geometry called also dual
pairs, play a relevant role in the theory of superintegrable
Hamiltonian systems. We prove the existence of an analogous
bif\/ibrated geometry in dynamical systems with a symmetry group
such that the reduced dynamics is periodic. The integrability of
such systems has been proven by M.~Field and J.~Hermans with a
reconstruction technique. We apply the result to the nonholonomic
system of a ball rolling on a surface of revolution.}


\Classification{37J35; 70H33}

\renewcommand{\thefootnote}{\arabic{footnote}}
\setcounter{footnote}{0}

\section{Introduction}

A number of nonholonomic mechanical systems are integrable, in the
sense that their motions are conjugate to a linear f\/low on the
torus, namely, are quasiperiodic. A noteworthy example is a heavy
sphere which rolls without sliding inside a convex surface of
revolution with vertical symmetry axis. This system has symmetry
group $\SOTRE\per S^1$ and the reduced system was known to be
integrable already to Routh \bib{routh}, who showed that the
equations of motion can be solved by quadratures depending on the
solution of a linear ODE. The quasiperiodicity of the unreduced
system was proven much more recently by J.~Hermans \bib{hermans95}
(see also \bib{zenkov}) by means of a remarkable reconstruction
argument.

In fact, it has been shown by M.~Field \bib{field} and
J.~Hermans~\bib{hermans95} that, if a dynamical system is
invariant under the (free) action of a compact Lie group and if
the reduced dynamics is periodic, then the reconstructed dynamics
is quasi-periodic. Specif\/ically, there is a f\/ibration of a
certain `regular' subset of the phase space by invariant tori of
dimension $r+1$, where $r$ is the rank of the group\footnote{Some
notions on Lie groups are collected in the Appendix.}, on which
the f\/low is conjugate to a linear f\/low. Moreover, the
frequencies of these motions depend only on the integrals of
motion of the reduced system. As proven by N.T.~Zung
in~\cite{zng06}, a stronger result is true in the Hamiltonian
setting, where quasi-periodicity of the reconstructed f\/low
follows also from quasi-periodicity of the reduced f\/low.

The purpose of this paper is to investigate the global geometry of
the f\/ibration by the inva\-riant tori of these systems, with the
analogy of the Hamiltonian case in mind.  Integrability of
Hamiltonian systems is usually related to the existence of
integrals of motion with certain proper\-ties via the
Liouville--Arnol'd theorem and its `noncommutative'
generalizations (see e.g.~\bib{karasev-maslov,fasso-cg} for
reviews and further references). The resulting structure is
characterized by the existence of not just one, but two invariant
foliations of the phase space, which are naturally def\/ined and
interrelated: the f\/ibration by the (isotropic) invariant tori
and its (coisotropic) `polar' foliation, which together form what
in symplectic geometry is called a `dual pair'. For our purposes,
it is especially useful to describe this structure in the case in
which the integrability is related to the existence of a symmetry
group (see e.g.~\bib{mischenko78,fomenko-trofimov,blaom} for
details, extensions, and further references).  Assume that the
Hamiltonian is invariant under the action of a compact and
connected Lie group~$G$, which acts freely and in a Hamiltonian
way on the phase space $M$. Assume moreover, just for simplicity,
that the momentum map $J:M\to\mathfrak{g}$ has connected f\/ibers.
If $\dim M=\dim G+\rank G$, then

\bList \bull The f\/ibers of the momentum map
$J:M\to\mathfrak{g}^*$ are the orbits of a maximal torus of $G$
and hence are dif\/feomorphic to $\toro^r$, where $r=\rank G$.
They are isotropic, and the dynamics is conjugate to a linear
f\/low on them. \bull Globally, the momentum map $J:M\to
B:=J(M)\subset \mathfrak{g}^*$ is a f\/ibration with f\/iber
$\toro^r$. \bull The orbits of $G$, being symplectically
orthogonal to the level sets of $J$, are coisotropic and union of
invariant tori.  More precisely, each $G$-orbit is a
$\toro^r$-bundle over a regular coadjoint orbit of $G$. \bull
Globally, the $G$-orbits are the f\/ibers of a f\/ibration
$\pi_A:M\to A$, where $A$ is a manifold of dimension equal to
$\rank G$ which can be indentif\/ied with a Weyl chamber of
$\mathfrak g$. The components of $\pi_A$ play the role of the
`actions' of the system. \bull The two f\/ibrations are
compatible, in the sense that $\pi_A=\rho\circ J$, where
$\rho:\mathfrak{g}^*\to A$ gives the foliation of $\mathfrak{g}^*$
into its coadjoint orbits. \eList

Sometimes, the Hamiltonian is independent of the momentum map and
the energy-momen\-tum map $(H,J)$ turns out to be the momentum map
of a $G\per S^1$-action (see noticeably \bib{blaom}, chapter~12).
In these cases, the previous picture applies with $G$ replaced by
$G\per S^1$ and $r$ replaced by $r+1$. This is the situation that
most directly compares to the case treated in the sequel.

The above `bif\/ibrated' structure, which plays an important role
in a number of questions, see e.g. \bib{karasev-maslov,fasso-cg},
can be described via the existence of a commutative diagram such
as in Fig.~1 left. Fig.~1 right gives a pictorial representation
of this structure, see e.g.~\bib{fasso-cg}, where the individual
f\/lowers are the coisotropic leaves, dif\/feomorphic to $G$, the
petals are the invariant tori $\toro^r$, the centers of the
f\/lowers are the coadjoint orbits of $\mathfrak{g}^*$, and the
meadow is a Weyl chamber, namely, the `action space'. The picture
is inaccurate in that it suggests that both f\/ibrations $J$ and
$\pi_A$ are topologically trivial, while it is not necessarily so.

\begin{figure}[ht]
\centerline{\includegraphics[width=10cm,height=3.5cm]{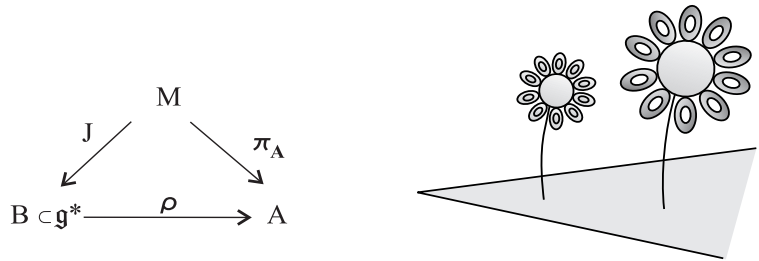}}
\caption{The bif\/ibration (left) and its pictorial representation
(right).} \label{fiori}
\end{figure}

In this paper, extending previous results by Hermans that we
review in Section 2, we show that a bif\/ibrated structure is
encountered also in the case of systems whose integrability is
determined via reconstruction from a periodic reduced dynamics
(Section~3). The dif\/ference, of course, is that there is no
symplectic interpretation of the `polarity' between the two
f\/ibrations. In Section~4 we relate the integrability to
Bogoyavlenskij's integrability~\bib{bog} and study the integrals
of motion produced by the reconstruction procedure. Finally, in
Section~5 we illustrate the previous constructions on the example,
already considered by Hermans, of the sphere inside a~convex
surface of revolution.

While preparing this paper, we were informed by Richard Cushman
that a very similar extension of Hermans' work will appear
in~\bib{CDS}.  The main dif\/ference, to our understanding, is
that we focus on the bif\/ibrated structure, while \bib{CDS}
emphazises non-compact groups and non-free actions.

\section{Reconstruction of the invariant tori}

{\bf A. The phase map.} For a symmetric system with periodic
reduced dynamics the reconstruction of the unreduced dynamics is a
well known topic, see e.g. \bib{MaMoRa}, and requires introduction
of a Lie-group valued map called {\it phase}. Consider: \vskip1mm
\bList{\it \item[{\rm H1.}] A free action $\Psi: G\per M\to M$ of
a compact connected Lie group $G$ on a manifold $M$. }\eList
\vskip1mm Since the group is compact, freeness of the action
implies properness, so the quotient space $\bar M:=M/G$ is a
manifold. Moreover, the quotient map $\pi:M\to\bar M$ def\/ines a
$G$-principal bundle.  We often write $\Psi_g(m)$ or $g.m$ for
$\Psi(g,m)$ with $g \in G$ and $m \in M$. Next, we consider
\vskip1mm \bList{\it \item[{\rm H2.}] A $G$-invariant vector field
$X$ on $M$. }\eList \vskip1mm If we denote by $\Phi^X_t$ the map
at time $t$ of the f\/low of $X$, then
$\Psi_g\circ\Phi^X_t=\Phi^X_t\circ\Psi_g$ for all $t \in \R$ and
$g \in G$. The f\/low of the reduced vector f\/ield $\bar
X:=\pi_*X$ on $\bar M$ is intertwined by $\pi$ with that of $X$,
that is, $\pi\circ \Phi^{X}_t= \Phi_t^{\bar X}\circ\pi$. We
further assume that: \vskip1mm \bList {\it \item[{\rm H3.}] The
reduced vector field $\bar X$ has periodic flow, with positive
smooth (minimal) period. }\eList \vskip1mm Here, smoothness of the
period means that the function $\bar \tau:\bar M\to\R$, which
associates to each point $\bar m\in\bar M$ the minimal period
$\bar \tau (\bar m)$ of the $\bar X$-orbit through it, is smooth.
The function $\bar\tau$ def\/ines a function $\tau = \bar\tau
\circ \pi$ on $M$.

Hypothesis H3 also implies that the orbit space of $\bar X$, that
we denote $A$, is a smooth manifold, and that the quotient map
$\bar\pi:\bar M\to A$ def\/ines a $S^1$-principal bundle, see
\bib{fgs05}. We denote $\pi_A$ the function $\bar\pi \circ
\pi:M\rightarrow A$, as in the following diagram:
$$
\xymatrix{
&& A\\
M\ar[r]^\pi\ar@/_/[rrd]_{\tau}\ar@/^/[rru]^{\pi_A}& \bar M \ar[ru]_{\bar\pi}\ar[rd]^{\bar\tau}\\
&& \R }
$$

The def\/inition of the phase map is based on the fact that the
periodicity of the reduced dynamics implies that the unreduced
orbit of a point $m$ in $M$ returns periodically, with
period~$\tau(m)$, to the $G$-orbit through $m$, see
Fig.~\ref{faseetori}.  Therefore, there exists an element
$\gamma(m) \in G$ (unique, since the action is free) such that
\begin{gather*}
  \P^X_{\tau(m)}(m) \ug \Psi_{\gamma(m)}(m) .
\end{gather*}
This def\/ines the {\it phase map} $\gamma:M\to G$. This map is
called monodromy in \bib{hermans95}.

In reconstruction theory, where $\gamma(m)$ is usually called
total phase or reconstruction phase, the emphasis is often on
computation and properties of the individual reconstructed orbits.
A~dif\/ferent approach has been pursued by M.~Field \bib{field}
and by J.~Hermans \bib{hermans95} who used this map, together with
the periodicity of the reduced f\/low, to prove the
quasi-periodicity of the unreduced dynamics. A central role in
this approach is played by the fact that the phase map is a
conjugacy-equivariant $G$-valued integral of motion of the
unreduced dynamics:

\begin{proposition}
Under hypotheses {\rm H1--H3}, the phase map $\gamma:M\to G$ is
\bList \item[\rm (i)]  smooth; \item[\rm (ii)] constant along the
$X$-orbits, that is
  $\gamma\circ\P^X_t=\gamma$ for all $t$;
\item[\rm (iii)] equivariant with respect to the group conjugacy,
that is
 $\gamma(g.m) = g \gamma(m) g^{-1}$ for all $g \in G$ and $m \in M$.
\eList
\end{proposition}

\begin{proof}
  Consider the manifold $\mathcal M=\{(m,g.m)\, |\, m\in M\quad g\in
  G\}\subset M\times M$ and the submersion $f:M\times G\rightarrow
  \mathcal M$ def\/ined by $(m,g)\mapsto (m,g.m)$. Since $\mathcal
  N=\{(m,\P^X_{\tau(m)}(m)\,|\,m\in M\}$ is a smooth submanifold of
  $\mathcal M$ of codimension $\dim G$, the implicit function theorem
  together with the fact that projection onto $M$ of $f^{-1}(\mathcal
  N)$ is $1:1$ implies that $f^{-1}(\mathcal N)$ is locally the graph of
  a smooth function $\gamma: M\rightarrow G$.

  Statements (ii) and (iii) follow from the facts that the
  period function $\tau$ is constant along the orbits of $X$
  and along the orbits of $G$, that is $\tau(m)=\tau(\P^X_t(m))=\tau(g.m)$ for all $g \in G$ and $t \in
  \R$. It follows that $\P^X_{\tau(m)}(\P^X_t(m)) =
  \P^X_t(\P^X_{\tau(m)}(m)) = \P^X_t(\gamma(m).m) =
  \gamma(m).(\P^X_t(m))$, hence $\gamma(\P^X_t(m))=\gamma(m)$,
  and that $\P^X_{\tau(m)}(g.m)= g. \P^X_{\tau(m)}(m)=
  g.\gamma(m).m=(g\gamma(m)g^{-1})(g.m)$, hence
  $\gamma(g.m)=g\gamma(m)g^{-1}$.
\end{proof}

{\bf B. Invariant tori.} Let $Z(h)$ be the centralizer of $h \in
G$. As is recalled in the Appendix, for a dense and open set of
points $h\in G$, called regular points, $Z(h)$ is a torus
contained in $G$ of maximal dimension $r$. We denote $M_\reg$ the
set of points $m\in M$ whose phase $\g(m)$ is a regular element of
$G$.

\begin{figure}[ht]
\centerline{\includegraphics[width=7cm]{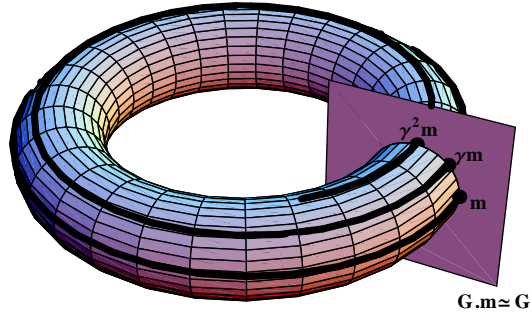}}
\caption{The $T_m$-orbit of the centralizers of $\g(m)$ (the
circle) in the $G$-orbit (the shaded plane), and its transport
along the f\/low of $X$ (the torus).} \label{faseetori}
\end{figure}

For $m\in M_{\rm reg}$, the maximal torus $T_m:=Z(\g(m))$, and
hence the set $T_m.m$, is dif\/feomorphic to $\toro^r=(\R/\Z)^r$.
Since all the elements of $T_m.m$ have the same phase, and since
$T_m.m$ is mapped onto itself by the map $\Phi^X_{\tau(m)}$, one
can def\/ine a foliation of $\Mreg$ by $X$-invariant tori of
dimension $r+1$ by transporting the $T_m$-orbits along the f\/low
of $X$, see Fig.~\ref{faseetori}.  More precisely, following
\bib{field,hermans95}, we have:

\begin{proposition}
Consider a point $m \in \Mreg$. Let $\eta_m$ be an element in  the
Lie algebra $\CartanAlg_m$ of $T_m$ such that $\g(m)\ug
\exp\eta_m$. Define the map
\begin{gather*}
  j_m : S^1\per T_m \;\to\; M_\reg , \qquad j_m(\a,h)
        = \Psi_{h} \circ \P^X_{\a\tau(m)} \circ
          \Psi_{\exp(\a\eta_m)}^{-1} (m) .
\end{gather*}
\bList \item[\rm (i)] The map $j_m$ is an embedding and its image
$P_m := j_m(S^1\per T_m)$ is invariant under the flow of $X$ and
diffeomorphic to $\toro^{r+1}$. \item[\rm (ii)] If $m'$ belongs to
$P_m$ then $P_{m'}=P_m$. \eList
\end{proposition}

\begin{proof} (i) Since $S^1\per T_m$ is compact we need only to show that $j_m$ is
  an injective immersion. To prove that the map is an immersion
  consider a vector $(a, \xi)\in T_{(\bar\alpha,\bar h)}(S^1\times
  T_m)$, with $a\in\R$ and $\xi\in \CartanAlg_m$. Such a vector can be
  represented by the curve $t\mapsto (\bar\alpha+ta,\exp(t\xi)\bar h)$
  whose image via $j_m$ is $t\mapsto \Psi_{ \exp(t\xi)\bar h} \circ
  \P^X_{(\bar\alpha + t a)\tau(m)} \circ \Psi_{\exp((\bar\alpha + t
  a)\eta)}^{-1} (m)$. By dif\/ferentiation with respect to $t$, one sees
  that $(a,\xi)$ is mapped onto $X_\xi+a\tau(m) X-a X_\eta$, with
  $X_\xi,X_\eta$ the inf\/initesimal group action associated to the
  elements $\xi,\eta$. This vector vanishes if and only if $a=0$ and
  $\xi=0$.

  The injectivity follows from the fact that if
  $j_m(\alpha,h)=j_m(\beta,k)$ then $\P^X_{(\a-\beta)\tau(m)}(m) =
  \Psi_{h^{-1}\exp\big((\alpha-\beta)\eta\big) k}(m)$. The left
  hand-side belongs to $G.m$ if and only if $n=\alpha-\beta\in\Z$, in
  which case the equation becomes $\gamma(m)^n.m =
  h^{-1}\gamma(m)^nk.m$. This equation holds if and only if $h=k$.

  The invariance of $P_m$ under the f\/low of $X$ is obvious. In fact if
  $m'=j_m(\alpha,h)$, then $t\mapsto j_m\Big(\alpha+t,
  \exp\big(t\,\eta(m)\big)h\Big)$ is the f\/low of $X$ through
  $m'$.

(ii) Let $m'=j_m(\alpha,h)\in P_m$ for some $\alpha\in S^1$ and
$h\in T_m$,
  hence $m=j_{m'}(-\alpha,h^{-1})$. If $m''\in P_m$, then
  $m''=j_m(\beta,k)$ and hence $m''=j_{m'}(kh^{-1},\beta-\alpha)\in P_{m'}$.
\end{proof}

We shall call the sets $P_m$ either `invariant tori of $X$' or,
more frequently, {\it petals}. Note that while, for given $m$, the
embedding $j_m$ depends on the choice of the element $\eta_m$ in
the algebra, the petal through it depends only on $X$ and on the
$G$-action.

Item ii. in the previous Proposition implies that the petals give
a partition of the manifold~$M_\reg$. However, since in every
suf\/f\/iciently small open set $U\subset M_\reg$ one can choose
$\eta_m$ as a~smooth function of $m$, the argument used in the
proof of the Proposition yields a $\toro^{r+1}$-action on $U$
whose orbits are the petals $P_m$, $m\in U$. This implies that the
partition of $M_\reg$ in petals is a foliation. We do not detail
this construction here since we shall later show that the petals
are in fact the f\/ibers of a (global) f\/ibration.

\medskip

{\bf C. Quasi-periodicity of the f\/low.} Following
\bib{hermans95}, we now show that the f\/low of $X$ on each petal $P_m$
is quasi-periodic, that is, it is conjugate to the linear f\/low
on $\toro^{r+1}$. To this end we introduce suitable coordinates on
the maximal torus $T_m$. Specif\/ically, choose a basis
$\xi_1,\ldots,\xi_r$ of the Lie algebra $\CartanAlg_m$ of $T_m$
that generates the lattice of elements of $\CartanAlg_m$ that
exponentiate to the identity. The map $\Xi :
(\beta_1,\ldots,\beta_r) \mapsto \exp\Big(\sum_j\beta_j\xi_j\Big)$
is a dif\/feomorphism from $\toro^r$ onto $T_m$ and
\begin{gather*}
  i_m \;:=\; j_m\circ(\id\times\Xi) \;:\;
             S^1\per \toro^r \;\to\; \Mreg
\end{gather*}
is an embedding of $\toro^{r+1}$ into $\Mreg$. Explicitly,
$i_m(\a,\beta)= \Psi_{\Xi(\beta)} \circ \P^X_{\a \tau(m)} \circ
\Psi_{\exp(\a\eta)}^{-1} (m)$.

\renewcommand{\mod}{\mathrm{(mod\;1)}}

\begin{proposition}
Consider $m\in\Mreg$ and let $\g(m)=\exp\eta$ with
$\eta=\sum_j\eta_j\xi_j$. The map $i_m$ intertwines the flow of
$X$ and the linear flow on~$\toro^{r+1}$ given by
\begin{gather*}
     \vp_t(\a,\beta) \;:=\;
     \left(\a+\frac{t}{\tau(m)} \,,\;
          \beta+\frac{t}{\tau(m)}\vec\eta\right)
     \quad \mod ,
\end{gather*}
where $\vec\eta=(\eta_1,\ldots,\eta_r)$.
\end{proposition}

\begin{proof} This is a simple computation, let $\tau=\tau(m)$:
\begin{gather*}
  i_m\big(\vp_t(\a,\beta)\big)
  =
  \Psi_{(\Xi(\beta+\frac t\tau\vec\eta))}
  \circ
  \P^X_{(\a+\frac t\tau)\tau}
  \circ
  \Psi^{-1}_{\exp((\a+\frac t\tau)\eta)}
  (m)
  \\
 \phantom{i_m\big(\vp_t(\a,\beta)\big)}{}=
  \Psi_{\Xi(\beta)}
  \circ \Psi_{\exp(\frac t\tau\eta)}
  \circ \P^X_{\a \tau}
  \circ \P^X_{t}
  \circ \Psi^{-1}_{\exp(\a\eta)}
  \circ \Psi_{\exp(\frac t\tau\eta)}^{-1}
  (m)
  \\
\phantom{i_m\big(\vp_t(\a,\beta)\big)}{}=
  \P^X_{t}
  \circ \Psi_{\Xi(\beta)}
  \circ \P^X_{\a \tau}
  \circ \Psi^{-1}_{\exp(\a\eta)}
  (m)
  \\
\phantom{i_m\big(\vp_t(\a,\beta)\big)}{}=
  \P_t^X(i_m(\a,\beta)),
\end{gather*}
where we have used the commutativity of $X$-f\/low and $T$-action
and the fact that the latter is Abelian.
\end{proof}

\section[The bifibration]{The bif\/ibration}

\noindent{\bf A. The f\/lower f\/ibration.} Since the dynamics in
two petals $P_m$ and $P_{g.m}=g.P_m$ is conjugate by~$g$, it is
quite natural to group the petals to form bigger sets of
$G$-related petals
\begin{gather*}
    F_m \,:=\, G.P_m ,\qquad m\in M_\mathrm{reg} ,
\end{gather*}
that we call {\it flowers}. A hint about the relevance of this
further structure can be found in \bib{hermans-thesis}, where
however it is not exploited.

We begin the study of this structure by characterizing the
`internal' structure of each f\/lower, that is, its quotient over
the petals. Note that, by Proposition 1, each f\/lower is
contained in~$\Mreg$ and the maximal tori of its points are all
conjugate.

\begin{proposition}
Consider $m\in M_\reg$ and let $T_m$ be the maximal torus
$Z(\gamma(m))$.  Then, $F_m$ is diffeomorphic to $S^1 \times G$
and is a $\toro^{r+1}$-principal bundle with the petals as fibers
and base $G/T_m$.
\end{proposition}

\begin{proof} The proof that $F_m \cong S^1\per G$ is very
similar to the proof of Proposition 2. Fix $m\in M_\reg$ and
choose $\eta_m$ such that $\exp(\eta_m)=\gamma(m)$. Consider the
map
\begin{gather*}
J_m:S^1 \times G\rightarrow M, \qquad (\alpha, g)\mapsto \Psi_g
\circ \P^X_{\a \tau(m)} \circ \Psi_{\exp(\a\eta_m)}^{-1} (m).
\end{gather*}
For the reasons given in Proposition 2 this map is an embedding,
so its image is dif\/feomorphic to $S^1\per G$. That this image is
$F_m$ follows from the fact that $F_m = \bigcup_{0\le
t<\tau(m)}\P_t^X(G.m)$.

The smooth map
\begin{gather*}
  \chi: S^1\per T_m \per F_m \to F_m, \qquad
  \chi(a,b , J_m\big(\alpha,g)\big) = J_m(a+\alpha, g b)
\end{gather*}
gives a free action of $\toro^{r+1}$ on $F_m$, which is therefore
a $\toro^{r+1}$-principal bundle. The orbits of this action are
the sets $g.P_m$, with $g\in G$. Observing that $J_m$ pulls back
each set $g.P_m$ to the set $S^1\times gT_m$ in $S^1\times G$, one
concludes that the quotient of $F_m$ by the petals is
dif\/feomorphic to the quotient of $S^1\times G$ by $S^1\times
T_m$.
\end{proof}

{\bf B. The petal-f\/lower bif\/ibration.} Proposition 4 implies
also that the f\/lowers are the f\/ibers of a f\/ibration of
$M_\mathrm{reg}$. In fact, they are the level sets of the map
\begin{gather*}
  \pi_A:M_\mathrm{reg} \to A_\mathrm{reg}:=\pi_A(M_\mathrm{reg})
  \subset A
\end{gather*}
and, by the Ehresmann f\/ibration theorem, any submersion with
compact and connected f\/ibers is a locally trivial f\/ibration,
see \bib{meigniez}. We now prove that also the petals are the
f\/ibers of a~f\/ibration of $M_\reg$, thus obtaining the
following bif\/ibrated structure:

\begin{theorem} There exists a manifold $B$ and a commutative diagram
$$
\xymatrix{
& M_\reg \ar[dl]_{\pi_B} \ar[dr]^{\pi_A} &  \\
B\ar[rr]^{\rho}& &A_\mathrm{reg}}\\
$$
where: \bList \bull the map $\pi_B$ is a fibration whose fibers
are the petals, diffeomorphic to $\toro^{r+1}$; \bull the level
sets of the fibration $\pi_A$ are the flowers, diffeomorphic to
$S^1\times G$; \bull the level sets of $\rho$ are the the quotient
space flowers/petals and are diffeomorphic to $G/T$, where $T$ is
any chosen maximal torus of $G$. \eList The dynamical system
evolves in the petals, and is there a linear flow.
\end{theorem}

The remainder of this section is devoted to the proof of this
Theorem, which is the main result of the paper. In
\bib{hermans95}, it is proven that the petals form the f\/ibers of
{\it semilocal} f\/ibrations, where semilocal means in a
neighbourhood of a petal.

\medskip

{\bf C. The semiglobal normal form. } The proof of Theorem 1 is
rather technical and is articu\-lated in two steps. First we
describe the {\it semiglobal} geometry of the f\/ibration by
petals, where semiglobal means in a $G$-invariant neighbourhood of
each f\/lower, not just in a neighbourhood of a petal. We do so by
def\/ining a sort of semiglobal normal form for the bif\/ibration
in petals and f\/lowers. The argument we use is a direct extension
of that of Hermans
\bib{hermans95}, who essentially proves the same result
semilocally. In the next subsection we shall construct the
base~$B$ of the f\/ibration by petals with a cut and paste
construction.

\begin{lemma}\label{Prop5}
For every $m\in M_\reg$ there is an open neighbourhood $U$ of
$\pi_A(m)$ in $A_\mathrm{reg}$, a discrete $\Z$-action on $G\per
U\per\reali$ and a diffeomorphism $\iota_U :M_U:=\pi_A^{-1}(U)\to
(G\times U\times \R)/\interi$ which intertwines: \bList \bull the
flow of $X$ in $M_U$ and the map $t'\mapsto [g,u,t+t']$ on
$(G\times U\times \R)/\interi$; \bull the action of $G$ on $M_U$
and the $G$-action $g'.[g,u,t]=[g'g,u,t]$ on $(G\times U\times
\R)/\interi$, \eList where the square brackets represent the
equivalence classes.
\end{lemma}

The proof of this Lemma rests upon another, technical, Lemma:

\begin{lemma}\label{sezionebuona}
  Given any maximal torus $T$ of $G$ and any $a\in A_\reg$, there
  exists an open neighbour\-hood $U$ of $a$ and a section
  $\sigma:U\rightarrow M_\reg$ of $\pi_A$ such that the phase map
  $\g\circ \sigma$ has values~in~$T$.
\end{lemma}

\begin{proof}[Proof of Lemma 2.] Fix $a\in A_\reg$. Since both maps
$\pi:M\rightarrow \bar M$ and $\bar\pi:\bar M\rightarrow A$ are
principal bundles there is a neighbourhood $U\subset A_\reg$ of
$a$ and a section $\sigma':U\rightarrow M_\reg$ of $\pi_A$. Let
$m=\sigma'(a)$. The element $\g(m)$ belongs to $G$ and necessarily
there exists an element $g\in G$ such that $\g(g.m)=g\g(m)g^{-1}$
belongs to $T$ (see the Appendix). It follows that the section
$\sigma''=g.\sigma'$ maps $a$ to the element $g.m$ such that the
phase $\g(g.m)\in T$.

  The Lie algebra $\mathfrak g$ of $G$ can be decomposed into
  $\CartanAlg \oplus V$, where $\CartanAlg$ is the Lie algebra of $T$ and
  $V$ is a subspace of $\mathfrak g$. Consider the map
\begin{gather*}
  f:V \times U \longrightarrow G,
  \qquad (\xi,u)\mapsto \Ad_{\exp(\xi)}\g(\sigma''(u)).
\end{gather*}
Observe that $f(0,a)$ belongs to $T$, that $T$ is a submanifold of
$G$ of codimension equal to the dimension of $V$, and that
$d_{(0,a)}f(V \times 0)$ is a subspace of the tangent space
$T_{\g(\sigma''(a))}G$ transversal to $T_{\g(\sigma''(a))}T$. The
implicit function theorem implies that possibly restricting $U$,
$f^{-1}(T)$ is the graph of a smooth function $\xi:U\rightarrow V$
such that $\xi(a)=0$.  This allows the def\/inition of a section
$\sigma:u\mapsto\exp(\xi(u))\sigma''(u)$, which has the required
property.
\end{proof}

\begin{proof}[Proof of Lemma 1.] Fix a maximal torus $T$ of $G$. There
is an open covering $\cU$ of $A_\reg$ by open sets $U$ as in Lemma
2, namely, for each $U\in\cU$ there is a section $\sigma_U$ of
$\pi_A$ such that $\g\circ\sigma_U(U)\subset T$.  Def\/ine the map
$\iota_U:G\times U\times \R\rightarrow M_\reg$, $(g,u,t)\mapsto
\Phi_t^X(g.\sigma_U(u))$. This map is easily shown to be an
immersion, invariant under the action of $G$, and invariant under
the f\/low of $X$. Consider an element $y$ in the image of
$\iota_U$ and a triplet $(g,u,t)$ such that $\iota_U(g,u,t)=y$.
The set of elements $(g',u',t')$ such that $\iota_U(g',u',t')=y$
is the set $\{(g\g(\sigma_U(u))^{-n},u,t+n \bar\tau(u))\,|\, n\in
\Z\}$. We hence def\/ine the $\Z$-action $n.(g,u,t) =
(g\g(\sigma_U(u))^{-n},u,t+n\bar\tau(u))$. The quotient of
$G\times U\times \R$ under this $\Z$-action is dif\/feomorphic to
the image of $\iota_U$, that is the open set $M_U=\pi_A^{-1}(U)$.

In each open set $M_U$ the map $\iota_U$ conjugates the f\/low of
$X$ to the $\R$-action $t'[g,s,t]=[g,s,t+t']$ and the $G$-action
to $g'[g,s,t]=[g'\, g, s, t]$.
\end{proof}

We call the dif\/feomorphism $\iota_U$ between $M_U$ and $(G\times
U\times \R)/\Z$ a \emph{semiglobal normal form} for $M$.

The proof of Lemma 1 shows also that the f\/low of $X$ is
contained in the sets $\iota_U(g T\times u\times \R)$, which are
the petals $P_{\iota_U(g,u,0)}$, and in the bigger sets $\iota_U(G
\times u\times \R)$, which are the f\/lowers $F_{\iota_U(g,u,0)}$.
We can thus write the semiglobal commutative diagrams
$$
\begin{array}{lr}
\xymatrix{
& M_U \ar[dl]_{\pi_{B_U}} \ar[dr]^{\pi_A}\\
G/T\times U\ar[rr]& &U} \qquad\qquad & \xymatrix{
& [g,u,t] \ar[dl]\ar[dr]\\
[g]_T\times u\ar[rr]& &u}\end{array}
$$
where $[g]_T$ stands for the class $g T$ in $G/T$. The f\/ibers of
$\pi_{B_U}$ are the petals, the f\/ibers of~$\pi_A$ are the
f\/lowers. The quotients obtained identifying the petals of a
f\/lower to points are all dif\/feomorphic to $G/T$. We call such
manifolds \emph{centers}.

\begin{proof}[Proof of Theorem 1.] In order to prove Theorem 1, we now
globalize the above construction of the maps
$\pi_{B_U}:M_U\rightarrow G/T\times U$ so as to def\/ine a
manifold $B$ and a map $\pi_B:M_\reg \rightarrow B$. The
construction of $B$ requires the preliminary construction of a
certain covering space $\widetilde A$ of $A_\reg$ by means of a
cut and paste construction that we now describe.

Let $\mathcal U$ be a covering of $A_\reg$ by open sets $U$ as in
Lemma 1 and let $W=\{[e]_T=[g_0]_T,[g_1]_T,\dots\}$ be the Weyl
group $N(T)/T$ (see the Appendix). An open covering of $\widetilde
A$ is made of the open sets $W\times U$ with $U\in\cU$. The gluing
map between $W\times U_1$ and $W\times U_2$ is given by imposing
that $([g_1]_T,u_1)\in W\times U_1$ equals $([g_2]_T,u_2)\in
W\times U_2$ if and only if $u_1=u_2$ and $P_{g_1.\sigma_1(u_1)}
=P_{g_2.\sigma_2(u_2)}$. In the expression above, for each
$i=1,2$, $\sigma_i$ is the section of $\pi_A$ over $U_i$ whose
existence is stated in Lemma 1 and $g_i$ is any representative in
$N(T)$ of the element $[g_i]_T\in N(T)/T$. Note that the choice of
a dif\/ferent representative $g_i b$, $b\in T$, has no ef\/fects
on the petals since $P_{g_i b.\sigma(u)} =P_{g_i.\sigma(u)}$. In
conclusion, every element of $\widetilde A$ can be non-uniquely
represented as a pair $([g_i]_T,u)$ with $[g_i]_T\in W$ and $u\in
U$. The topology of $\widetilde A$ is induced by the topology of
the sets $W\times U$, the space is Hausdorf\/f and paracompact
because $A$ is Hausdorf\/f and paracompact and $W$ is f\/inite.

The Weyl group $W$ acts as deck transformations of $\widetilde A$.
Given an element of $\widetilde A$ represented by $([g_i]_T,u)$ in
$W\times U$, the element $[g_j]_T\in W$ maps it to
$([g_jg_i]_T,u)$. The Weyl group also acts on the space $G/T$ by
right-multiplication, $[g_j]_T [g]_T=[gg_j^{-1}]_T$. We can hence
def\/ine the manifold $B=G/T\times_W \widetilde A$, where the
f\/ibered product is obtained by means of the anti-diagonal action
$[g_j]_T.\Big([g]_T,([g_i]_T,u)\Big)=\Big([gg_j^{-1}]_T,([g_jg_i]_T,u)\Big)$.
The space $B$ is a covering space of $G/N(T)\times A$ with f\/iber
$W$ and is in turn covered by $G/T\times \widetilde A$ with
f\/iber $W$. Its points can be non-uniquely represented by pairs
$([g]_T,([g_i]_T,u))$.

We now def\/ine the map $\pi_B:M_\reg \rightarrow B$. For each
$m\in M_\reg$ there exists an open set $U\in\cU$ and three
elements $u\in U$, $g\in G$ and $t\in \R$ such that
$m=\iota_U(g,u,t)$, see Lemma~1.  We thus def\/ine
\begin{gather*}
  \pi_B(m) := \big([g]_T,([e]_T,u)\big) \in G/T\times_W \widetilde A,
\end{gather*}
where $e$ is the identity of $G$. In order to prove that the map
$\pi_B$ is well def\/ined we need to show that, if the point $m$
is also represented by $\iota_{U'}(g',u',t')$ with $u'\in
U'\in\cU$, then the element $([g']_T,([e]_T,u'))$ coincides with
$([g]_T,([e]_T,u))$. Both sets of data satisfy
\begin{gather}\label{fondamentale}
m=\Phi_{t'}^X\left(g'.\sigma'(u')\right)=\Phi_{t}^X\left(g.\sigma(u)\right),
\end{gather}
where $\sigma$ and $\sigma'$ are the sections of Lemma 1. This
implies that
$\sigma'(u')=\Phi_{t-t'}^X\left(g'^{-1}g.\sigma(u)\right)$. This
means that $([g'^{-1} g]_T,u)$ and $([e]_T,u')$ represent the same
point of $\widetilde A$. Hence $([g']_T, ([e]_T,u'))=([g']_T,
([g'^{-1} g]_T,u))$. By using the def\/inition of f\/ibered
product we f\/inally conclude that
\begin{gather*}
([g']_T, ([g'^{-1}
g]_T,u))=([g']_T[g^{-1}g']_T^{-1},[g^{-1}g']_T([g'^{-1}
g]_T,u))=([g]_T,([e]_T,u)).\tag*{\qed}
\end{gather*}\renewcommand{\qed}{}
\end{proof}

\section{The phase map and the integrals of motion}

As pointed out in particular by Bogoyavlenskij \bib{bog}, but see
also
\bib{fedorov}, quasi-periodicity of a f\/low can be (semilocally)
linked to the presence of a number of f\/irst integrals and of a
complementary number of commuting dynamical symmetries which
preserve these f\/irst integrals. Specif\/ically, if on a manifold
$M$ and for some $n<d=\dim M$ there are \bList \bull a submersion
$F:M\to\reali^{d-n}$ with compact and connected f\/ibers, and
\bull $n$ everywhere linearly independent and pairwise commuting
vector f\/ields $Y_1,\ldots,Y_n$ which are tangent to the f\/ibers
of $F$, \eList then the f\/ibers of $F$ are dif\/feomorphic to
$\tenne$ and any vector f\/ield $X$ on $M$ such that
\begin{gather*}
   L_XF_i=0  \qquad\mathrm{and}\qquad [X,Y_j]=0 , \qquad
   i=1,\ldots,d-n,\quad j=1,\ldots,n ,
\end{gather*}
is conjugate to a constant vector f\/ield on $\tenne$.

In the setting of the previous sections, where $n=r+1$, the
f\/ibration $\pi_B:M_\reg\to B$ ensures the existence of local
sets of $\dim M-r-1$ functionally independent integrals of motion
of~$X$, which are obtained by lifting local coordinates on $B$.
For what concerns the commuting dynamical symmetries, the vector
f\/ield $X$ and the chosen maximal torus $T$ def\/ine a semiglobal
group action of $T^{r+1}$ (see the proof of Proposition~4). This
action can be non-global, due to the possible non-trivial
monodromy of the torus bundle (see
\bib{duist80,cushman-duist} for the notion of monodromy). Nonetheless, the
semiglobal torus action def\/ines in every set $M_U$ a family of
$r$ independent vector f\/ields.

Going beyond these elementary considerations, we link now the
integrals of motion to the pha\-se map $\gamma$, which is a Lie
group-valued conserved quantity. Clearly, the projection $\pi_A:M
\rightarrow A$ gives rise to the semiglobal existence of sets of
$\dim A$ functionally independent integrals of motion of $X$,
which are obtained by lifting functions on $A$. We thus focus on
the integrals of motion coming from $B$ but not from $A$.

In the semiglobal normal form of Lemma 1, the expression of all
these integrals is immediate. Since the petals (which are the
level sets of a full set of integrals of motion) correspond to
subsets of the form $[gT,u,\R]$, the projection onto $u$ (the map
$\pi_A$) gives the integrals obtained by lifting functions on $A$.
The other integrals of motion can then be obtained via the map
\begin{equation}
  \delta_U: M_U \longrightarrow G/T, \qquad \[g,u,t\] \mapsto \[g\].
\end{equation}
In fact, two points $[g,u,t]$ and $[g',u',t']$ in $M_U$ belong to
the same f\/lower if and only if \mbox{$u'=u$}. They belong to the
same petal if and only if $u'=u$ and $[g]_T=[g']_T$. Globalizing
this construction, we have

\begin{proposition}
Fix any maximal torus $T$ in $G$ and let $N(T)$ be the normalizer
of $T$ in~$G$. Then, the integrals of motion that are not lifts of
functions on $A_\reg$ via $\pi_A$ are lifts of functions on
$G/N(T)$ via the map
\begin{gather*}
  \delta: M_\reg\rightarrow G/N(T) ,\qquad m \mapsto [g_m] ,
\end{gather*}
where $g_m\in G$ is such that $g_m^{-1}\gamma(m)g_m\in T$.  The
level sets of this map consist of $|N(T)/T|$ petals.
\end{proposition}

\begin{proof} Given $U\in\cU$, take any $m\in \pi_A^{-1}(U)\subset
  M_\reg$ and let $g_m$ be such that $g_m^{-1}\gamma(m)g_m$ belongs to
  $T$. Then, there exists an element $h\in N(T)$ and a time $t$ such
  that $\Phi^X_t(hg_m^{-1}.m)=\sigma_U(\pi_A(m))$. It follows that the
  point $m$ is represented in a semiglobal normal form of Lemma 1 by
  the triple $[g_m h^{-1}, \pi_A(m), -t]$. Hence, $\delta(m)= [g_m
  h^{-1}]_T$ in $G/T$. A further projection on $G/N(T)$ gives the
  element $[g_m]_{N(T)}$.
\end{proof}

\begin{remark} The function $\delta$ described in Proposition 5
is obtained by composing the function~$\delta_U$ in formula (2)
with the projection $G/T\rightarrow G/N(T)$. This gives a map that
has the same rank as~$\delta$. The only dif\/ference is in the
number of connected components of the level sets. This choice
allows an easier def\/inition of~$\delta$.
\end{remark}

\begin{remark} The manifold $G/T$ is dif\/feomorphic to any regular coadjoint
orbit, while $G/N(T)$ is dif\/feomorphic to a coadjoint orbit
modulo the action of the Weyl group. This gives a very strict
analogy with the case of a Hamiltonian system with a Hamiltonian
group action we described in the Introduction.
\end{remark}

\section{An example: the ball rolling on a surface of revolution}

We illustrate now the above construction in the example of a heavy
sphere rolling without sliding inside a convex surface of
revolution which has vertical axis and faces upward. This system
has been considered by Hermans, who actually used it as motivation
and exemplif\/ication of his reconstruction theory
\bib{hermans95}. We add a global perspective to this.
For general reference on this system see
\bib{hermans95,zenkov,fgs05} and references therein.

The phase space of the system is the eight-dimensional manifold
$\R^2\per\R^2\times \SO(3) \times \R \ni (a,\dot a,Q,w)$, where
$a\in\rdue$ are Cartesian coordinates of the center of mass of the
sphere, $Q\in\SOTRE$ f\/ixes the attitude of the sphere, and
$w\in\R$ is the component of the angular velocity normal to the
surface. The dynamical vector f\/ield of the system is invariant
under the action of $G=S^1\times\SO(3)$ given by $(\vt,R).(a,\dot
a,Q,w)=(S_\vt\, a,S_\vt\,\dot a, RQ,w)$, where $S_\vt$ is the
rotation of angle $\vt$ about the vertical axis. This action is
free on the submanifold $M=\left(\R^4\setminus\{0\}\right) \times
\SO(3) \times \R$ and the four-dimensional quotient manifold $M/G$
is dif\/feomorphic to $\bar M = \left(R^3 \setminus \{0\}\right)
\times \R$, with projection $(a,\dot a,w)\mapsto (b(a,\dot a),w)$,
where $b:\reali^4\setminus\{0\}\to\R^3\setminus\{0\}$ is the Hopf
f\/ibration.

It was classically known that this system has three independent
integrals of motion, the energy~$H$ and two functions $J_1$ and
$J_2$ which are sometimes referred to as Chaplygin integrals.
These three integrals are $G$-invariant and hence descend to
independent integrals of motion $\bar H$, $\bar J_1$, $\bar J_2$
of the reduced system $\bar X$ in the quotient manifold $\bar M$.
Moreover, the level sets of $\bar H$, $\bar J_1$, $\bar J_2$ are
compact, and hence are f\/inite disjoint union of circles, see
\bib{zenkov,fgs05}. These circles are the orbits of the reduced
dynamical system $\bar X$, which is therefore periodic. (This was
also shown using a dif\/ferent approach in \bib{hermans95}). The
continuity of the period is shown in
\bib{fgs05}. Thus, the system f\/its in the setting of the
previous sections and, since $S^1\per\SOTRE$ has rank two, the
dynamics is quasi-periodic on three-dimensional tori, with
frequencies which depend only on $H$, $J_1$, $J_2$
\bib{hermans95}. Thus, there are two additional integrals of
motions, besides $H$, $J_1$, $J_2$.

A complete description of these integrals, and of the
bif\/ibration, faces a number of dif\/f\/iculties. First,
determination of the phase map appears to be prohibitive even for
given (and simple) prof\/iles of the surface of revolution.
Without it, it is obviously not possible to completely determine
$M_\reg$, the two base manifolds $B$ and $A_\reg$, and the two
additional integrals of motion as in Proposition 5. Moreover,
since it is also dif\/f\/icult to establish whether the f\/ibers
of $(H,J_1,J_2)$ are connected, we do not know if $A=\pi_A(M)$
coincides with the image of $M$ under the map $(H,J_1,J_2)$ or if
it is a covering space of this image. Nonetheless, something can
be said about the structure of the bif\/ibration in this case:

\begin{proposition}
  The regular region $M_\reg$ of the phase space of the ball rolling on a
  surface of revolution admits a petal-flower bifibration in which:
\bList \bull the petals are $3$-tori \bib{hermans95}; \bull the
five-dimensional base $B$ is a $2$-sheeted covering of the
manifold $\R P^2\times A_\reg$; \bull the flowers are
diffeomorphic to $\SO(3) \times \tdue$ and the fibration in petals
of each flower is a $\toro^3$-bundle over a $2$-sphere,
$\SO(3)\times\tdue \rightarrow S^2$, with projection the Hopf map
from the $\SO(3)$-component onto $S^2$. \eList
\end{proposition}

\begin{proof}
The only fact we need to prove is the statement about the
structure of $B$. Since the Weyl group is in this case $\Z_2$,
$B=S^2\times_{\Z_2} \widetilde A_\reg$ with $\widetilde A_\reg$ a
$2$-sheeted covering of $A_\reg$. Thus, $B$~is a two sheeted
covering of $\R P^2\times A_\reg$.
\end{proof}

Once we have partly clarif\/ied the structure of $B$, we may
describe, up to the determination of $\g$ and some f\/inite
covering, the two additional functions of $B$ that are not lifts
of functions of~$A_\reg$. To this ef\/fect we need to def\/ine two
independent functions on $\R P^2$, lift them to $B$, and then lift
them again to $M_\reg$. Note that the full set of lifted functions
on $M_\reg$ might have some extra f\/inite symmetry, in which case
their levels will be f\/inite families of petals.

We need to re-trace what was done in Proposition 5 in order to
def\/ine the function $\delta$. Any point $(a,\dot a, Q,w)$ in $M$
has phase $\gamma(a,\dot a, Q,w)=(\vt(a,\dot a, Q,w),R(a,\dot a,
Q,w))$, with $\vt$ a circle-valued function and $R$ an
$SO(3)$-valued function. Following Proposition 5, we f\/ix a
maximal torus of $S^1\times SO(3)$: we elect the torus $S^1\times
S^1_z$, where $S^1_z$ are the rotations around the $z$-axes. Thus,
$\delta$~is a map from $M_\reg$ into $S^1\per\SOTRE/N(S^1\per
S^1_z) \cong\R P^2$.  The element $R(a,\dot a, Q,w)$ is a
rotation~$R^\vp_{\vec v}$ of angle $\vp$ about the oriented axis
$\vec v$. A rotation that conjugates $R^\vp_{\vec v}$ to
$R^\vp_{\vec e_3}$ is $R_{\vec v\times \vec e_3}^{\arccos(\vec
v,\vec e_3)}$. It follows that $\g(R_{\vec v\times \vec
e_3}^{\arccos(\vec v,\vec e_3)}.(a,\dot a, Q, w)) \in S^1\times
S^1_z$, and hence belongs to the maximal torus we have chosen: we
have determined the element $g_m$ of Proposition 5. We can hence
state:

\begin{proposition}
The map $\delta$ has values in $\R P^2$ and is given by $(a,\dot
a, Q,w)\mapsto [R_{\vec v\times\vec e_3}^{\arccos(\vec v,\vec
e_3)}\vec e_3]$. (In this last expression, the square brackets
denote the equivalence classes with respect to the projection
\hbox{$\R^3\setminus\{0\} \to \R P^2$}.)
\end{proposition}

\section{Appendix: used facts on compact Lie groups}

Recall that the {\it centralizer} of an element $h\in G$ is the
subset of $G$ of the elements that leave $h$ f\/ixed,
\begin{gather*}
   Z(h)= \{g\in G \,:\; ghg^{-1}=h\}.
\end{gather*}
It is immediate to prove that the centralizer is a subgroup. The
{\it normalizer} of a subgroup $H$ of $G$ is the set of elements
that leave $H$ f\/ixed
\begin{gather*}
  N(H)=\{g\in G \,:\; gHg^{-1}=H\}.
\end{gather*}
$N(H)$ is a subgroup of $G$ which contains $H$ and is the biggest
subgroup of $G$ in which $H$ is a~normal subgroup.

Given a compact group $G$, and $h\in G$, the subgroup $Z(h)$ is a
compact subgroup of $G$ and, if Abelian, is a maximal torus. When
this happens, $h$ is called {\it regular}, $Z(h)\cong T$ is called
a~\emph{Cartan subgroup}, and its Lie algebra $\CartanAlg$ is
called a \emph{Cartan subalgebra}. The dimension $r$ of a~Cartan
subgroup is the same for all Cartan subgroups and is called the
\emph{rank} of the group $G$. All Cartan subgroups are conjugate
(see page~159 in \bib{dieckbrocker}).

Consider now a Cartan subgroup $T$ and its Cartan subalgebra
$\CartanAlg$. The normalizer of $T$ in $G$ is a f\/inite extension
of $T$ (see page~158 in \bib{dieckbrocker}), and the f\/inite
group $W=N(T)/T$ is called the {\it Weyl group}.

The group $N(T)$ acts by conjugation on $T$, action that obviously
def\/ines an action of the Weyl group on $T$ and, by
linearization, an action on the Cartan algebra $\CartanAlg$. The
Weyl group action on $\CartanAlg$ is a group of ref\/lections (see
page~192 in
\bib{dieckbrocker}), each ref\/lection in the Weyl group f\/ixes an
hyperplane, called a {\it Weyl wall}.

This fact, together with the fact that every element of $G$
belongs to at least one Cartan subgroup, implies that the set of
regular elements is an open dense set. To complete the picture we
observe the trivial fact that the exponential of elements in a
Weyl wall give precisely the elements which belong to more than a
maximal torus or, in other words, those whose centralizer is a
bigger, non-commutative subgroup of $G$.

\subsection*{Acknowledgements}

The authors thank the Bernoulli Center (EPFL, Lausanne) for its
hospitality during the 2004 program {\em Geometric Mechanics and
Its Applications}, where the biggest part of this work was done,
and Hans Duistermaat for some enlightening conversations on these
topics.

\pdfbookmark[1]{References}{ref}
\LastPageEnding

\end{document}